\documentclass[12pt]{article}
\usepackage{amsfonts,amssymb} 
\usepackage{amsmath}
\numberwithin{equation}{section}
\begin{document}
\title{A unified approach to the integrals of Mellin--Barnes--Hecke type}
\author{
Gopala Krishna Srinivasan
\footnote{e-mail address: gopal@math.iitb.ac.in}
}
\date{
 Department of  Mathematics, Indian Institute of Technology Bombay, 
Mumbai 400 076
}
\maketitle
\begin{center}
{\bf Abstract} 
\end{center}
In this paper we provide a unified approach to a family of integrals of Mellin--Barnes type using distribution 
theory and Fourier transforms.   
 Interesting features arise in many of the cases  which call for the 
application of pull-backs of distributions via smooth submersive maps defined by H\"ormander. 
We derive by this method  the 
integrals of Hecke and Sonine relating to various types of Bessel functions which have found 
applications in analytic and algebraic number theory. 
\begin{center}
\end{center}
\section{Introduction:} The theory of distributions and Fourier transforms    
has been successfully applied in the 
theory of differential operators to obtain precise asymptotic properties of solutions. 
In this paper we look at some 
applications of distribution theory in the study of special functions by providing a unified approach to a certain 
class of integrals of a type studied by Mellin and Barnes.

Integrals involving products of gamma functions along vertical lines were studied first by Pincherle in 1888 
and an extensive theory was developed by Barnes \cite{barnes} and Mellin \cite{mellin}. 
Cahen \cite{cahen} employed some of these integrals  
in the study of the Riemann zeta function and other Dirichlet series.  
 In the spirit of Mellin's 
theory some of Ramanujan's formulas have been generalized by G. H. Hardy \cite{hardy}, p. 98 ff. 
The work of Pincherle provided impetus for the subsequent investigations of Mellin \cite{mellin} and Barnes \cite{barnes} on the integral representations of 
solutions of generalized hypergeometric series (see \cite{nielsen}, chapter 16 and 
the comment on p. 225).   
A detailed commentary 
on Pincherle's work \cite{pincherle} set against a historical backdrop is available in \cite{mainardi}.

Among the integrals studied by Barnes, the integral formula \eqref{barnes_first_integral} 
is well known which served as the point of departure for 
Barnes for his development of the theory of the hypergeometric functions in his seminal paper \cite{barnes}.  
A more exotic example \eqref{barnes_second}-\eqref{barnes_second_integral} 
 appeared in a later paper by Barnes.    
The integrals of Mellin--Barnes also play an important role in the theory of the q-analogues of the hypergeometric functions introduced by Heine in 1847.   
For a discussion of Mellin--Barnes integrals we refer to \cite{whittaker1}, p. 286 ff and the thorough investigations  in  
the recent books \cite{andrews} and \cite{paris}.  
 For the classical evaluation of these and
 other integrals of this class see \cite{andrews} particularly pp. 89-91 and pp. 151-154. 
A complete account of q-hypergeometric series and their
 Mellin--Barnes integral representations is available in \cite{gasper}, chapter 4.

In his work on real quadratic number fields (see \cite{hecke}, p. 349)   
Hecke employs the transformation formula \eqref{hecke} of a Mellin--Barnes integral which was generalized by Rademacher with a 
view towards applications to number fields of higher degree  
(see \cite{rademacher}, p. 58).  
The integral of Hecke also features in his proof of Hamburger's theorem on the 
Riemann zeta function (see \cite{hecke}, p. 378) where its relation to the  Bessel function of the 
third kind, the Hankel functions, is established.   
Hecke's formula is reminiscent of an integral (\eqref{sonine} below) considered by Sonine in his researches on the Bessel
 function. 
The formulas of Sonine and Hecke combine to yield an integral representation of the Macdonald function $K_p(x)$ 
that is occasionally employed in analytic number theory.

We provide a transparent and unified approach  to these results in the present paper using 
 Fourier integrals with non-linear phase functions that yield distributions (densities) 
given by the pull-back of the Dirac delta distribution. 
We present in sections 3 and 4 proofs of all the formulas stated in the introduction (besides a few others)  
in the spirit of Fourier analysis. It turns out that all these formulas fall out of 
the basic equation proved in section 2:
\begin{equation}\label{funda-formula}
\int_{-\infty}^{\infty}dt\int_{\mathbb R^k}\exp{(-it\rho(u_1, u_2,\dots, u_k))}\phi(u_1,\dots,u_k)\;du_1\dots du_k = 
2\pi\langle\rho^{*}(\delta_0), \phi\rangle.             
\end{equation}
Here $\rho(u_1, u_2,\dots, u_k)$ is a phase function, $\phi$ is an arbitrary 
member of the Schwartz class ${\cal S}(\mathbb R^k)$ 
of rapidly decreasing functions (see \cite{hormander}, chapter VII), 
$\delta_0$ denotes the standard Dirac delta distribution and $\rho^{*}(\delta_0)$ 
denotes the 
pull-back of Dirac delta distribution by $\rho$ (see \cite{hormander}, p. 136 or
 \cite{duistermaat}, p. 103). 
It is useful to write \eqref{funda-formula} in the more suggestive notation 
\begin{equation}
\int_{-\infty}^{\infty}\exp{(-it\rho(u_1, u_2,\dots, u_k))}\;dt = 2\pi\delta(\rho(u_1,\dots, u_k) = 0). 
\end{equation}
A noteworthy special case of (1.2) is the \it Fourier--Gel'fand formula \rm \cite{duistermaat}, p. 193: 
\begin{equation}\label{fourier-gelfand}
\int_{\mathbb R}\exp{it(u-v)}\;dt = 2\pi\delta(u-v). 
\end{equation}

The paper is organized as follows. Section 2 contains the proof of the basic formula \eqref{funda-formula} 
that we repeatedly use. 
Section 3 contains the proofs of the formulas of Ramanujan and their generalizations due to Barnes, Mellin and Hardy. 
Section 4 is devoted to 
the proof of  Hecke's formula and its relation to various types of Bessel functions. 
Finally, a few other examples susceptible to the same treatment, 
are gathered up in theorem (4.7). 
\section{A basic formula from distribution theory:}  
We shall prove in this section the requisite formula from the theory of distributions that we shall 
frequently employ. 
We recall first the notion of pull-back of distributions via smooth maps defined in \cite{hormander}, pp 134-136.  
For general background on distribution theory we also recommend the recent comprehensive work \cite{duistermaat}.   
\paragraph*{Theorem 2.1 (Pull-back of distributions by smooth submersions):} 
\begin{itemize}
\item[(i)]
Let $X$ and $Y$ be open sets in $\mathbb R^n$ and $\mathbb R^m$ respectively and  
$\rho : X \longrightarrow Y$ be a smooth submersive map. Then there exists a unique continuous linear map 
$\rho^{*}: {\cal D}^{\prime}(Y) \longrightarrow {\cal D}^{\prime}(X)$ such that $\rho^{*}(u) = u\circ\rho$ for all $u \in C(Y)$. 
\item[(ii)]
In the special case when $X = \mathbb R^n$ and $Y = \mathbb R$ the pull-back $\rho^{*}(\delta_0)$ is given by  
\begin{equation}
\rho^{*}(\delta_0) = dS/|\nabla \rho|              
\end{equation}
where $dS$ is the Euclidean surface measure along the zero locus $\rho(x_1, x_2,\dots, x_n) = 0$. 
\end{itemize}

With this notion we state and prove the main formula we use. 
\paragraph*{Theorem 2.2:} 
Suppose that $\rho : \mathbb R^n \longrightarrow \mathbb R$ 
is a smooth submersive map and $\phi(x)$ is a member of the Schwartz 
class ${\cal S}(\mathbb R^n)$,  then
\begin{equation}\label{funda-formula2}
\int_{-\infty}^{\infty} dt \int_{\mathbb R^n}\phi(x)\exp(it\rho(x)) \;dx = 2\pi\langle \rho^{*}(\delta_0), \phi
\rangle                                                                 
\end{equation}
\paragraph*{Proof:} Assume first that $\phi$ has compact support along which $\frac{\partial\rho}{\partial x_1} > 0$.  
Denoting the integral in \eqref{funda-formula2} by $I$, 
\begin{eqnarray*}
I &=& \lim_{\epsilon\rightarrow 0}\int_{-\infty}^{\infty} dt \int_{\mathbb R^n}\phi(x)\exp(it\rho(x) - \epsilon t^2)\;dx\\
&=& \lim_{\epsilon\rightarrow 0}\int_{\mathbb R^n}\phi(x) dx \int_{-\infty}^{\infty}\exp(it\rho(x) - \epsilon t^2)\;dt\\
&=& \sqrt{\pi}\lim_{\epsilon\rightarrow 0}\frac{1}{\sqrt{\epsilon}}\int_{\mathbb R^n}\phi(x)\exp(-(\rho(x))^2/4\epsilon)\;dx\\
\end{eqnarray*}
The assumptions on $\rho$ justify the change of variables\footnotemark\footnotetext{
The implicit function theorem gives only a neighborhood on which this holds but we may begin by assuming
 that the support of $\phi$ is contained in this neighborhood and then use a partition of unity.}
$$
\Psi \;:\; y_1 = \rho(x_1,x_2,\dots, x_n), y_2 = x_2, \dots, y_n =  x_n,   
$$ 
and the first of these defines a function $x_1 = \psi(y_1, x_2, \dots, x_n)$. The  
integral now transforms into 
$$
I = \sqrt{\pi}\lim_{\epsilon\rightarrow 0}\frac{1}{\sqrt{\epsilon}}\int_{\mathbb R^n}
\phi(\psi(y_1, y_2, \dots, y_n), y_2,\dots, y_n) \Big(\frac{\partial\rho}{\partial x_1}
\Big)^{-1}e^{-y_1^2/(4\epsilon)}\;dy_1dy_2\dots dy_n,
$$
where $\frac{\partial\rho}{\partial x_1}$ is evaluated at $(\psi(y_1, y_2, \dots, y_n), y_2,\dots, y_n)$. Rescaling the 
variable $y_1$,  
\begin{eqnarray*}
I = 2\sqrt{\pi}\lim_{\epsilon\rightarrow 0}\int_{\mathbb R^n}
\frac{\phi(\psi(2\sqrt{\epsilon} y_1, y_2, \dots, y_n), y_2,\dots, y_n)}{\partial_{1}\rho(\psi(2\sqrt{\epsilon}y_1, y_2,\dots, y_n), y_2,\dots, y_n)} 
e^{-y_1^2}\;dy_1dy_2\dots dy_n  \\
= 2\pi\int_{\mathbb R^{n-1}}\frac{\phi(\psi(0, y_2,\dots, y_n), y_2,\dots, y_n)}{\partial_{1}\rho(\psi(0, y_2,\dots, y_n), y_2,\dots, y_n)} 
\;dy_2dy_3\dots dy_n,
\end{eqnarray*}
by the dominated convergence theorem and smoothness of the functions involved. 
Now the Euclidean area measure on the locus $\rho(y_1, y_2,\dots, y_n) = 0$ is given by 
$$
dS = \Big\{1 + \sum_{j=2}^{n}\Big(\frac{\partial \psi}{\partial y_j}\Big)^2\Big\}^{1/2}dy_2\dots dy_n = 
\frac{|\nabla\rho|}{\partial_1\rho}\;dy_2\dots dy_n,
$$
and invoking the preceding theorem we get 
$$
I = 2\pi \int_{\mathbb R^{n-1}}\phi(\psi(0, y_2,\dots, y_n), y_2,\dots, y_n) \frac{dS}{|\nabla\rho|} = 
2\pi \langle \rho^{*}\delta_0, \phi
\rangle.
$$
The assumption that $\frac{\partial\rho}{\partial x_1} > 0$ along the support of $\phi$ is easily removed by employing a 
  partition of unity. Finally if $\phi$ does not have compact support one can use a sequence $(\phi_m)$ of smooth 
 functions with compact supports converging to $\phi$ in ${\cal S}(\mathbb R^n)$ and then 
\begin{eqnarray*}
\int_{-\infty}^{\infty} dt \int_{\mathbb R^n}\phi(x)\exp(it\rho(x))\; dx &=& 
\lim_{m\rightarrow\infty}\int_{-\infty}^{\infty}dt\int_{\mathbb R^n}\phi_m(x)\exp(it\rho(x))\; dx \\
&=& \lim_{m\rightarrow\infty}2\pi\langle
\rho^{*}(\delta_0), \phi_m\rangle = 2\pi\langle \rho^{*}(\delta_0), \phi\rangle.\\
\end{eqnarray*}
\paragraph*{Remark:} The author is grateful to one of the referees for suggesting the following alternative 
proof of theorem (2.2) using some of the results in \cite{duistermaat} and \cite{duistermaat1}. The notation 
$\rho_*$ denotes the pushforward (see \cite{duistermaat}, p. 93 ff.) 
under the map $\rho$. 
From problem (10.24) and the Fourier inversion formula (14.31)
 in \cite{duistermaat} it follows that 
$$
\frac{2\pi}{\|\mbox{grad}\; \rho\|} \delta_{\rho^{-1}(\{0\})} = 2\pi \rho^{*}\delta_0 = 
\rho^{*}({\cal F}{\bf 1}_{\mathbb R}),
$$
where ${\bf 1}_{\mathbb R}$ denotes the characteristic function of $\mathbb R$ and ${\cal F}$ denotes the 
Fourier transform. Furthermore we have, for $\phi \in C_0^{\infty}(\mathbb R)$,
\begin{eqnarray*}
\langle \;\rho^{*}({\cal F}{\bf 1}_{\mathbb R}),\; \phi \rangle &=& 
\langle \;{\cal F}{\bf 1}_{\mathbb R},\; \rho_*\phi \rangle = 
\langle \;{\bf 1}_{\mathbb R},\; {\cal F}(\rho_*\phi) \rangle = \int_{\mathbb R}{\cal F}(\rho_*\phi)(\xi) d\xi \\
\end{eqnarray*} 
Using the explicit formula (10.23) in \cite{duistermaat} for the pushforward,  
\begin{eqnarray*}
\langle 2\pi\rho^{*}\delta_0, \phi\rangle &=& 
\int_{\mathbb R}\int_{\mathbb R} e^{-i\xi x}
\int_{\rho^{-1}(\{x\})}\frac{\phi(y)\; dy \;dx\; d\xi}{\|\mbox{grad}\;\rho(y)\|} \\
&=& \int_{\mathbb R}\int_{\mathbb R} \int_{\rho^{-1}(\{x\})} \frac{e^{-i\xi\rho(y)}\;\phi(y)\; dy \;dx\; d\xi}{\|\mbox{grad}\;\rho(y)\|}\\
\end{eqnarray*}
Using now the result of exercise (7.36) in \cite{duistermaat1}, the last integral transforms into 
$$
\int_{\mathbb R}\int_{\mathbb R^n}e^{-i\xi\rho(x)}\phi(x)dx\;d\xi. 
$$ 
and we get the desired formula (2.2). 
\section{Integrals of Barnes, Mellin, Ramanujan and Hardy:} 
For details on the gamma function we refer to \cite{duistermaat} (pp. 164-168), \cite{nielsen}, \cite{srinivasan} and 
\cite{whittaker1} (chapter 12) and record here for convenience the basic results that we shall use. 
\paragraph*{Theorem 3.1:} (i) Duplication Formula: If Re$\;a > 0$, 
$\displaystyle{
\;\sqrt{\pi}\;\Gamma(a) =  2^{a-1}\Gamma\Big(\frac{a}{2}\Big)\Gamma\Big(\frac{a+1}{2}\Big).
}$ 
\vskip 0.08 in 
(ii) Euler's reflection formula:
$\displaystyle{
\Gamma(a)\Gamma(1-a) = \frac{\pi}{\sin\pi a}, \quad a \notin \mathbb Z. 
}$
\vskip 0.1 in
(iii) For $c > 0$ fixed, the function $t \mapsto |\Gamma(c + it)|$ decays exponentially fast as $t\rightarrow \pm \infty$ 
namely,  
\begin{equation}\label{stirling}
|\Gamma(c + it)| = |t|^{c-\frac{1}{2}}\exp(- \pi |t|/2  + O(1)), \quad t \rightarrow \pm \infty
\end{equation}
\paragraph*{Proof:} We shall briefly indicate a proof of (iii) using the Stirling's formula 
(see \cite{whittaker1}, p. 248 or \cite{binet} pp. 244-245): 
$$
\log\Gamma(c+it) = \big(c+it -\frac{1}{2}\big)\log (c+it) - it + O(1),\quad |t|\rightarrow \infty. 
$$
Taking the real part of the above equation and exponentiating we get 
$$
|\Gamma(c + it)| = (c^2+t^2)^{\frac{c}{2}-\frac{1}{4}}\exp(-t\tan^{-1}(t/c) + O(1)),
\quad |t|\rightarrow \infty
$$
from which \eqref{stirling} follows. A different proof is given on 
pp. 238-239 of \cite{pincherle} (see also p. 234).

We recall here the definition of the beta function:
\begin{equation}\label{beta}
B(p, q) = \int_0^1 u^{p-1}(1-u)^{q-1} du,\quad \mbox{Re}\;p > 0,\;\;\mbox{Re}\;q > 0.
\end{equation}
The basic relation between the beta and gamma function is given by 
\begin{equation}\label{beta-gamma}
\Gamma(p+q)B(p, q) = \Gamma(p)\Gamma(q), \quad \mbox{Re}\;p > 0,\;\;\mbox{Re}\;q > 0.
\end{equation}
We begin with a result due to Binet (\cite{binet}, p. 136) that we shall frequently use. 
\paragraph*{Theorem 3.2 (Binet):} When Re$\;p > 0$ and Re$\; q > 0$, 
\begin{equation}\label{binet}
B(p, q) = \int_{-\infty}^{\infty}\frac{(e^{(p-q)u} + e^{(q-p)u})\;du}{(e^u + e^{-u})^{p+q}}.  
\end{equation}
\paragraph*{Proof:} The result is obtained by substituting $u = e^{x/2}(e^{x/2}+e^{-x/2})^{-1}$ 
in the integral \eqref{beta} and using the symmetry $B(p, q) = B(q, p)$.

Ramanujan gave the following formula for the Fourier transform of $|\Gamma(a+it)|^2$ for positive 
real values of $a$ which we now state and prove using (1.1) (see also \cite{debraj}).
\paragraph*{Theorem 3.3 (Ramanujan):} If $a$ is a complex number in the open right half plane, then the following holds. 
\begin{equation}\label{ramanujan1}
\int_{-\infty}^{\infty} \Gamma(a + it)\Gamma(a - it) e^{-i\xi t} dt = 
{\sqrt{\pi}}\Gamma(a)\Gamma(a + \frac{1}{2})\mbox{sech}^{2a}\xi/2, 
\end{equation}
\paragraph*{Proof:} 
Using \eqref{beta}, \eqref{beta-gamma} and \eqref{funda-formula2} in succession we get    
\begin{equation*}
\int_{-\infty}^{\infty}\Gamma(a +it)\Gamma(a-it)e^{-i\xi t} dt  
=  
2\Gamma(2a)\int_{-\infty}^{\infty} dt\int_{-\infty}^{\infty}\frac{e^{it(2u-\xi)} du}{(e^u + e^{-u})^{2a}}  
= 4\pi\Gamma(2a)\langle\rho^{*}\delta_0, (e^u + e^{-u})^{2a}\rangle,
\end{equation*}
where $\rho(u) = 2u - \xi$ and so $\rho^{*}(\delta_0) = \frac{1}{2}\delta_{\xi/2}$. Thus the integral in 
\eqref{ramanujan1} equals 
$$
\frac{2\pi\Gamma(2a)}{(4^a\cosh^{2a}\xi/2)}.
$$  
Applying the duplication formula we get the desired result.  
\paragraph*{Theorem 3.4 (Barnes' first integral):} For $a, b, c$ and $d$ in the right half plane we have 
\begin{equation}\label{barnes_first_integral}
\frac{1}{2\pi}\int_{-\infty}^{\infty}\Gamma(a+it)\Gamma(b+it)\Gamma(c-it)\Gamma(d-it)\;dt = 
\frac{\Gamma(a+c)\Gamma(a+d)\Gamma(b+c)\Gamma(b+d)}{\Gamma(a+b+c+d)}.        
\end{equation}
\paragraph*{Proof:} Assuming $a, b, c$ and $d$ are real and positive, we  
 split the four gamma factors into groups of two and obtain using Binet's formula,  
\begin{eqnarray*}
\Gamma(a+it)\Gamma(c-it) &=& 
\Gamma(a+c)\int_{-\infty}^{\infty}\frac{\cosh(a-c +2it)u \;du}{(e^u + e^{-u})^{a+c}} \\
\Gamma(b+it)\Gamma(d-it) &=& \Gamma(b+d)\int_{-\infty}^{\infty}\frac{\cosh(b-d +2it)v \;dv}{(e^v + e^{-v})^{b+d}}. \\
\end{eqnarray*}
From this  we arrive at the following 
\begin{eqnarray*}
\int_{-\infty}^{\infty}\Gamma(a+it)\Gamma(b+it)\Gamma(c-it)\Gamma(d-it) \;dt = \\
\Gamma(a+c)\Gamma(b+d)
\int\int_{\mathbb R^2}\frac{dudv}{(2\cosh u)^{a+c}(2\cosh v)^{b+d}}
\int_{-\infty}^{\infty}E(t) \;dt.
\end{eqnarray*}
Here $E(t)$ is the sum of four terms $e(t, u, v) + e(t, u, -v) + e(t, -u, v) + e(t, -u, -v)$ and $e(t)$ is give by 
$\displaystyle{
e(t) = \exp({(a-c)u + (b-d)v +it(2u + 2v)}). 
}$  Using \eqref{funda-formula} we conclude 
\begin{eqnarray*}
\int_{-\infty}^{\infty} E(t)\;dt = e^{(a-c)u + (b-d)v}\delta(2u+2v = 0) + e^{(a-c)u - (b-d)v}\delta(2u-2v = 0) +\\
e^{-(a-c)u + (b-d)v}\delta(2u-2v=0) + e^{-(a-c)u - (b-d)v}\delta(2u+2v = 0).
\end{eqnarray*}
So we have finally the integral 
\begin{eqnarray*}
\int_{-\infty}^{\infty}\!\!\Gamma(a+it)\Gamma(b+it)\Gamma(c-it)\Gamma(d-it) \;dt &=& 
2\pi\Gamma(a+c)\Gamma(b+d)\int_{-\infty}^{\infty}\!\frac{\cosh(a-c-b+d)u du}{(2\cosh u)^{a+b+c+d}}\\
&=& \frac{2\pi\Gamma(a+c)\Gamma(b+d)\Gamma(a+d)\Gamma(b+c)}{\Gamma(a+b+c+d)}\\
\end{eqnarray*}
As a noteworthy special case we obtain by virtue of the duplication formula the following result \cite{ramanujan1} and 
\cite{hardy}, p. 103: 
\paragraph*{Corollary 3.5 (Ramanujan):} If $a$ and $b$ are positive real numbers, 
\begin{equation}\label{ramanujan2}
\frac{1}{\sqrt{\pi}}\int_{-\infty}^{\infty}|\Gamma(a+it)\Gamma(b+it)|^2 \;dt = 
\frac{\Gamma(a)\Gamma(a+\frac{1}{2})\Gamma(b)\Gamma(b+\frac{1}{2})\Gamma(a+b)}{\Gamma(a+b+\frac{1}{2})} 
\end{equation}
We now turn to an example where the gamma factor appears in the denominator (see \cite{hardy}, p. 103). 
\paragraph*{Theorem 3.6 (Ramanujan):} Suppose Re$\;(b-a) > \frac{1}{2}$, then 
\begin{equation}\label{ramanujan3}
\frac{1}{\sqrt{\pi}}\int_{-\infty}^{\infty}
\frac{\Gamma(a+it)\Gamma(a-it)}{\Gamma(b+it)\Gamma(b-it)}\;dt = 
\frac{\Gamma(a)\Gamma(a+\frac{1}{2})\Gamma(b-a-\frac{1}{2})}{\Gamma(b)\Gamma(b-\frac{1}{2})\Gamma(b-a))}.
\end{equation}
\paragraph*{Proof:} 
To first show that the integrand in \eqref{ramanujan3} is in 
$L^1(\mathbb R)$. Using \eqref{stirling},
$$
\Big|\frac{\Gamma(a+it)\Gamma(a-it)}{\Gamma(b+it)\Gamma(b-it)}\Big| \leq C|t|^{\mbox{2\small Re}(a-b)}\exp\Big(t\tan^{-1}(t/b) 
-t\tan^{-1}(t/a) 
\Big) = O(|t|^{\mbox{2\small Re}(a-b)}),
$$
from which the assertion follows. We now write 
$$
\frac{\Gamma(a+it)}{\Gamma(b+it)} = \frac{B(a+it, b-a)}{\Gamma(b-a)},\quad
\frac{\Gamma(a-it)}{\Gamma(b-it)} = \frac{B(a-it, b-a)}{\Gamma(b-a)},\quad 
$$ 
and hence the required integral is given by 
$$
\frac{4}{(\Gamma(b-a))^2}
\int\int_{\mathbb R^2}\frac{dudv}{(e^u + e^{-u})^{b}(e^v + e^{-v})^b}
\int_{-\infty}^{\infty}\frac{(\cosh(2a-b+it)u)(\cosh((2a-b-it)v)}{(e^u+e^{-u})^{it}(e^v+e^{-v})^{-it} }\;dt 
$$
The integral with respect to $t$ may be written as  
$$
\int_{-\infty}^{\infty}(E(t, u, v) + E(t, -u, v) + E(t, u, -v) + E(t, -u, -v))\;dt,
$$
where 
$E(t, u, v) = e^{(2a-b)u+(2a-b)v}\exp\big\{it(u-v-\log(e^u+e^{-u})+\log(e^v+e^{-v})\big\}.$ 
Note that the phase functions are not linear in $u$ and $v$. The zero locus of the 
phase function $\rho(u, v)$ appearing in $E(t, u, v)$ namely,  
$$
\rho(u, v) = u-v-\log(e^u+e^{-u})+\log(e^v+e^{-v}),
$$ 
is the line $u = v$ along which 
$\displaystyle{
|\nabla\rho(u, v)| = 2\sqrt{2}e^{-u}/(e^u + e^{-u})}$ whereby  
we get for the pull-back the weighted Lebesgue measure 
$\displaystyle{
\rho^{*}(\delta_0) = e^{u}(e^u + e^{-u})du/2. 
}$ 
The contribution to the integral from $E_1(t)$ is then 
$$
\frac{\pi}{(\Gamma(b-a))^2}
\int_{\mathbb R}\frac{e^{(4a-2b)u}e^{u}\;du}{(e^u + e^{-u})^{2b-1}}
$$
Likewise we determine the contributions from $E_2(t), E_3(t), E_4(t)$ and upon adding we get the value of 
Ramanujan's integral \eqref{ramanujan3} as 
$$
\frac{2\pi}{(\Gamma(b-a))^2}\int_{-\infty}^{\infty} \frac{(e^{(4a-2b+1)u} + e^{-(4a-2b+1)u})\;du}{(e^u + e^{-u})^{2b-1}}.
$$ 
Using Binet's expression \eqref{binet} and the 
duplication formula we get \eqref{ramanujan3}.

A more general result may be proved along the same lines:
\paragraph*{Theorem 3.7:} If Re$\;(b-a) \geq 0$, Re$\;(d-c) \geq 0$, Re$\;(b+d) > 1$ and Re$\;(b+d-a-c) > 1$ then, 
\begin{equation}
\frac{1}{2\pi}\int_{-\infty}^{\infty} \frac{\Gamma(a+it)\Gamma(c-it)}{\Gamma(b+it)\Gamma(d-it)}\;dt = 
\frac{\Gamma(a+c)\Gamma(b+d-a-c-1)}{\Gamma(b+d-1)\Gamma(b-a)\Gamma(d-c)}              
\end{equation}
We now state the Fourier transform of $\Gamma(a+it)/\Gamma(b+it)$ from which we can get another proof of the previous result 
by appealing to Parseval's formula.
\paragraph*{Theorem 3.8:} If Re$\;a > 0$ and Re$\;(b-a) > 0$ then $\Gamma(a+it)/\Gamma(b+it)\in L^p(\mathbb R)$ where 
$\{\mbox{Re}\;(b-a)\}^{-1} < p \leq \infty$. 
The Fourier transform is the function in $L^q$ where $q^{-1} > 1-\mbox{Re}\;(b-a)$ 
and is given by   
\begin{equation}\label{hardy1}
\frac{1}{2\pi}\int_{-\infty}^{\infty}\frac{\Gamma(a+it)}{\Gamma(b+it)}e^{-i\xi t} \;dt 
= \left\{\begin{array}{lll}
e^{\xi a}(1 - e^{\xi})^{b-a-1}(\Gamma(b-a))^{-1},& & \xi < 0\\
0. & & \xi > 0\\
\end{array} \right. 
\end{equation}
\paragraph*{Remark:} Integration by parts reveals that, for $\xi \neq 0$, the integral in \eqref{hardy1}
 is  convergent (though not always absolutely) using 
the following result (see \cite{whittaker1}, p. 251): 
$$
\frac{\Gamma^{\prime}(z)}{\Gamma(z)} = \log z - \frac{1}{2z} -2\int_0^{\infty}\frac{udu}{(u^2+z^2)(e^{2\pi u}-1)},
\quad \mbox{Re}\;z > 0. 
$$
Integral \eqref{hardy1} is equivalent to formula \eqref{hardy3} obtained by  G. H. Hardy (see \cite{hardy}, p. 98) 
 which we discuss next. Integral \eqref{hardy3} generalizes a formula of Ramanujan. 
\paragraph*{Theorem 3.9 (Hardy):} (i)   If $-p < a < q-p$, the following holds:
\begin{equation}\label{hardy2}
\frac{1}{2\pi}\int_{-\infty}^{\infty}y^{-(a+it)}\Gamma(a+p+it)\Gamma(q-p-a-it) \;dt = 
\frac{\Gamma(q)y^p}{(1+y)^q} 
\end{equation}
(ii) If $a > 0 $ and $0 < x < 1$ then, 
\begin{equation}\label{hardy3}
\frac{1}{2\pi}\int_{-\infty}^{\infty} \frac{\Gamma(a+it)x^{-a-it}\;dt}{\Gamma(a+q+it)} = (1-x)^{q-1}/\Gamma(q)  
\end{equation}
If $x > 1$ the value of the integral is zero. 
\paragraph*{Proof:} (i)   We readily transform the integral in \eqref{hardy2} in the form (1.2):   
$$
\frac{y^{-a}}{2\pi}\int_0^{\infty}\int_0^{\infty} e^{-u-v}u^{p+a-1}v^{q-p-a-1}\;dudv\int_{-\infty}^{\infty} 
\exp(it\rho(u, v))\;dt.
$$
The phase $\rho(u, v) = \log u - \log v - \log y$ vanishes along the line $u = vy$ 
 whereby $\rho^{*}(\delta_0) = vydv$ and the integral is 
$$
y^p \int_0^{\infty} e^{-v(1+y)} v^q \;dv = y^p\Gamma(q)/(1+y)^{q}
$$
as desired. 

(ii) Expressing the integrand in \eqref{hardy3} in terms of the beta function and using Binet's formula \eqref{binet},    
the integral transforms to 
\begin{eqnarray*}
\frac{x^{-a}}{2\pi\Gamma(q)}\int_{-\infty}^{\infty}\frac{du}{(e^u+e^{-u})^{a+q}}
\int_{-\infty}^{\infty}\exp(-it(\log x + \log(e^u+e^{-u})) 
\Big[
e^{(a-q)u+itu} + e^{(q-a)u-itu}
\Big]\;dt\\
\end{eqnarray*} 
Denoting by $\rho(u)$ either of the phase functions $\log x \mp u + \log(e^u+e^{-u})$,  we have
$$
|\nabla\rho(u)| = 2|1-x|
$$ 
along the locus $\rho(u) = 0$ which is simply the point given by $\exp(\pm u) = x/(1-x)$. 
With the Dirac measure concentrated at this point with 
normalization $(1-x)/2$ we get for the integral the value $(1-x)^{q-1}/\Gamma(q)$ as asserted. When $x > 1$ the locus $\rho(u) = 0$ is empty which 
explains why the integral vanishes in this case.
\paragraph*{Mellin--Barnes integrals for the Kummer functions:} 
The Kummer functions $\phantom{}_{1}F_1[a, b, x]$ and $U(a, b, x)$ are normalized linearly independent 
solutions of the confluent hypergeometric differential equation and have the following integral 
representation (see \cite{slater} pp 34-38): 
\begin{eqnarray}
\phantom{X}\phantom{}_{1}F_1[a, b, x] &=& 
\frac{\Gamma(b)}{\Gamma(b-a)\Gamma(a)}\int_0^1e^{xt}t^{a-1}(1-t)^{b-a-1}dt,\quad \mbox{Re}\;b > \mbox{Re}\;a > 0,
\label{eqn: kummer1} \\
\phantom{X}U(a, b, x) &=& \frac{1}{\Gamma(a)}\int_0^{\infty}e^{-xt}t^{a-1}(1+t)^{b-a-1}dt,\quad\mbox{Re}\;a > 0,\; \mbox{Re}\; x > 0.
\label{eqn: kummer2} 
\end{eqnarray}
\paragraph*{Theorem 3.10:} For real positive values of $x$ we have the following Mellin--Barnes integral 
representation for $\phantom{}_1F_1[a, b, x]$ and $U(a, b, x)$
\begin{eqnarray}
\phantom{}_{1}F_1[a, b, -x] &=& \frac{1}{2\pi}
\int_{-\infty}^{\infty}\frac{\Gamma(b)\Gamma(1-it)\Gamma(a-1+it)}{\Gamma(a)\Gamma(b-1+it)}x^{-1+it}\;dt,
\quad \mbox{Re}\;b > \mbox{Re}\;a.\label{eqn: kummer-mellin-barnes1}\\
U(a, b, x) &=& 
\frac{x^{1-a}}{2\pi}\int_{-\infty}^{\infty}\!\!
\frac{\Gamma(1-it)\Gamma(a-1+it)\Gamma(a-b+it)}{\Gamma(a)\Gamma(1+a-b)}x^{-it}\;dt, 
\label{eqn: kummer-mellin-barnes2}
\end{eqnarray}
where Re$\;(a-b+1) > 0$ in \eqref{eqn: kummer-mellin-barnes2}. 
\paragraph*{Proof:} Expressing the ratio of the gamma factors in the the first integral and using Binet's formula we get
$$
\frac{\Gamma(b)}{2\pi\Gamma(a)\Gamma(b-a)x}\int_{0}^{\infty}e^{-u}du\int_{-\infty}^{\infty}\frac{dv}{(e^v + e^{-v})^{b-1}}
\int_{-\infty}^{\infty}(E(t, u, v) + E(t, u, -v))\;dt,
$$
where 
$\displaystyle{
E(t, u, v) = \exp(it(\log x - \log u - \log(e^v+e^{-v})+v)\exp(2a-b-1)v.
}$ 
The zero locus of the phase function $\rho(u, v)$ is the curve $x = u(1+e^{-2v})$ along which $0 < u < x$ and 
$$
\frac{dS}{|\nabla\rho|} = \frac{1}{2}(1+e^{2v})du.
$$ 
The sum of the two integrals is then 
$$
\frac{x^{1-b}\Gamma(b)}{\Gamma(a)\Gamma(b-a)}\int_0^x e^{-u}u^{a-1}(x-u)^{b-a-1}du
$$ 
from which the result follows. For the case of $U(a, b, x)$, we get 
$$
U(a, b, x) = 
\frac{\Gamma(1+a-b)}{2\pi\Gamma(a)}\int_{0}^{\infty}e^{-u}u^{a-2}du\int_{-\infty}^{\infty}\frac{dv}{(e^v+e^{-v})^{a-b+1}}
\int_{-\infty}^{\infty} (E(t, u, v)+E(t, u, -v))\;dt, 
$$ 
where 
$\displaystyle{
E(t, u, v) = \exp((a-b-1)v)\exp(-it(\log x - \log u - 2v)). 
}$  
The zero locus of the phase function $\rho(u, v) = \log x - \log u -2v$ is the curve $x = ue^{2v}$ 
along which $dS/|\nabla\rho| = \frac{1}{2}du$. From this we immediately get the desired result.   
\paragraph*{Mellin--Barnes integrals for $F(a, b, c, z)$:} Introducing one more gamma factor in  
\eqref{eqn: kummer-mellin-barnes1} 
leads to the following Mellin--Barnes representation for the hypergeometric function 
(see \cite{barnes}, p. 142 or \cite{whittaker1}, p. 286).  
\paragraph*{Theorem 3.11:} For $z$ not on the negative real axis and Min$\;\{\mbox{Re}\;a, \mbox{Re}\;b\;\} > 1$, 
the following holds:
\begin{equation}\label{mellin-barnes-hypergeometric}
F(a, b, c, -z) = 
\frac{\Gamma(c)}{2\pi\Gamma(a)\Gamma(b)}
\int_{-\infty}^{\infty}\frac{\Gamma(a-1+it)\Gamma(b-1+it)\Gamma(1-it)z^{-1+it}\;dt}{\Gamma(c-1+it)}
\end{equation}
\paragraph*{Proof:} Taking $z = x$ to be real positive and proceeding as in the case of $\phantom{}_1F_1[a, b, x]$,   
the integral in  \eqref{mellin-barnes-hypergeometric} equals 
$$
\frac{1}{2\pi\Gamma(b)}\int_{0}^{\infty}e^{-u}du\int_0^{\infty}\frac{e^{-w}w^{b-2}dw}{B(a, c-a)}\!
\int_{-\infty}^{\infty}\frac{dv}{(e^v+e^{-v})^{c-1}}\int_{-\infty}^{\infty}(E(t, u, v, w) + E(t, u, -v, w)) \;dt,
$$
where 
$$
E(t, u, v, w) = \exp(it(\log(xw) - \log u + v - \log(e^v + e^{-v}))e^{(2a-c-1)v}.
$$  
The phase function $\rho(u, v, w)$ is  $\log(xw) - \log u + v - \log(e^v + e^{-v})$  and its zero locus is given by 
$xw = u(1+e^{-2v})$ which is non-empty only when 
$xw \geq u$. Expressing $v$ in terms of $u$ and $w$ we get 
$$
\frac{dS}{|\nabla\rho|} = \frac{xw dudw}{2(xw-u)},
$$ 
and the integral after a scaling transformation assumes the form 
\begin{eqnarray*}
\frac{1}{\Gamma(b)B(a, c-a)}\int_0^1 s^{a-1}(1-s)^{c-a-1}ds\int_0^{\infty}e^{-w(1+xs)}w^{b-1}dw = \\
\frac{1}{B(a, c-a)}\int_0^1 s^{a-1}(1-s)^{c-a-1}(1+xs)^{-b} ds.\\
\end{eqnarray*}
Using a familiar result (see \cite{lebedev}, p. 245), the last integral equals $F(a, b, c, -x)$. The condition that $z$ 
be real positive may be relaxed by appealing to analytic continuation. 
\paragraph*{Theorem 3.12 (Barnes' second integral):} Assume that $a, b, c, \lambda, \mu, \nu$ have  positive real parts and  satisfy $\nu = a + b + c + \lambda + \mu$. Then the 
integral 
\begin{equation}\label{barnes_second}
\frac{1}{2\pi}
\int_{-\infty}^{\infty}
\frac{\Gamma(a+it)\Gamma(b+it)\Gamma(c+it)\Gamma(\lambda-it)\Gamma(\mu-it)\;dt}{\Gamma(\nu+it)} 
\end{equation}
equals
\begin{equation}\label{barnes_second_integral}
\frac{\Gamma(\lambda+a)\Gamma(\lambda+b)\Gamma(\lambda+c)\Gamma(\mu+a)\Gamma(\mu+b)\Gamma(\mu+c)}{\Gamma(\nu-a)\Gamma(\nu-b)\Gamma(\nu-c)}
\end{equation}
\paragraph*{Proof:} 
We derive the result by assuming that $a, b, c, \lambda, \mu, \nu$ are real, positive and  $\mu > \lambda$. 
These may be relaxed later by appealing to analytic continuation. 
Proceeding as before, 
\begin{eqnarray*}
\frac{\Gamma(a+it)}{\Gamma(\nu+it)} &=& 
\frac{1}{\Gamma(\nu-a)}\int_{-\infty}^{\infty}\Big(
\frac{e^{(2a-\nu+it)v} + e^{-(2a-\nu+it)v}}{(e^v+e^{-v})^{\nu + it}}
\Big)dv  \\
\Gamma(b+it)\Gamma(c+it) &=& \int_0^{\infty}\int_0^{\infty} e^{-w_1-w_2}w_1^{b-1+it}w_2^{c-1+it}\;dw_1dw_2\\
\Gamma(\lambda-it)\Gamma(\mu-it) &=& \int_0^{\infty}\int_0^{\infty} e^{-u_1-u_2}u_1^{\lambda-1-it}u_2^{\mu-1-it}\;du_1du_2\\
\end{eqnarray*}
Thus \eqref{barnes_second} leads to an integral of the form (1.1) with phase function $\rho$ given by 
$$
\rho(v, u_1, u_2, w_1, w_2) = \log(w_1w_2) - \log(u_1u_2) \pm v - \log(e^v+e^{-v}),\quad u_1, u_2, w_1, w_2 > 0. 
$$
The zero locus $\rho(v, u_1, u_2, w_1, w_2) = 0$  is empty if $w_1w_2 \leq u_1u_2$. Along this zero-locus 
$$
\frac{dS}{|\nabla\rho|} = \frac{1}{2}(1+e^{2v})du_1du_2dw_1dw_2.
$$ 
By virtue of (1.1), the integral \eqref{barnes_second} equals 
$$
\frac{1}{\Gamma(\nu-a)}\int_{R}e^{-w_1-w_2-u_1-u_2}w_1^{b-\nu}w_2^{c-\nu}u_1^{a+\lambda-1}u_2^{a+\mu-1}(w_1w_2-u_1u_2)^{\nu-a-1}\;dw_1dw_2du_1du_2,
$$
where $R$ is the region in the positive orthant given by $R = w_1w_2 - u_1u_2 > 0$. After expressing it in the form 
$$
\int_0^{\infty} du_1\int_0^{\infty}dw_1\int_0^{\infty}dw_2
\int_0^{w_1w_2/u_1}\phi(u_1, u_2, w_1, w_2) \;du_2
$$
the integral is simplified in five steps. 
\begin{itemize}
\item[(i)] Perform a scaling transformation $u_2 = w_1w_2z/u_1$ with $z$ as the new variable ($0 < z < 1$). 
\item[(ii)] Integrate with respect to $w_2$.
\item[(iii)] Perform a scaling transformation $w_1 = su_1/z$ with $s$ as the new variable.
\item[(iv)] Integrate with respect to $u_1$.  
\item[(v)] Perform the transformation $s = (1-x)/x$ with $x$ as the new variable ($0 < x < 1$). 
\end{itemize}
We then get the following expression
$$
\frac{\Gamma(b+\lambda)\Gamma(c+\mu)}{\Gamma(\nu-a)}
\int_0^1z^{a+\lambda-1}(1-z)^{\nu-a-1}dz\int_0^1x^{\lambda+c-1}(1-x)^{b+\mu-1}(1-x(1-z))^{-b-\lambda}dx.
$$
The condition $\nu = a+b+c+\lambda+\mu$ and a familiar result on the hypergeometric function 
(see \cite{lebedev}, p. 245) yields
$$
\frac{\Gamma(b+\lambda)\Gamma(c+\lambda)\Gamma(b+\mu)\Gamma(c+\mu)}{(\Gamma(\nu-a))^2}
\int_0^1z^{a+\lambda-1}(1-z)^{\nu-a-1}F(b+\lambda, c+\lambda, \nu-a, 1-z)\;dz.
$$
The last integral is again expressible in terms of 
the hypergeometric function\footnotemark\footnotetext{
To apply the formula on p. 277 of \cite{lebedev} we have assumed $\mu > \lambda$. 
} (\cite{lebedev}, p. 277):
$$
\Big(\frac{\Gamma(a+\lambda)\Gamma(b+\lambda)\Gamma(c+\lambda)\Gamma(b+\mu)\Gamma(c+\mu)}{\Gamma(\nu-a)\Gamma(\lambda+\mu)}\Big)
F(b+\lambda, c+\lambda, \nu +\lambda, 1)
$$
On appealing to the well known Gauss--Euler formula for $F(a, b, c, 1)$  we get \eqref{barnes_second_integral}. 
\section{Integrals of Hecke and Sonine:} 
A very interesting integral transformation was discovered by Hecke in connection with 
certain problems in algebraic number theory and generalized by Rademacher (see \cite{rademacher}, p. 58). 
The formula was later used by Hecke in 
his proof of Hamburger's theorem on the characterization of the zeta function through the functional relation.  
We discuss in this section this integral and 
a closely related integral representation due to Sonine for the Bessel functions of the first kind. 
We also deduce from Sonine's and Hecke's formulas integral representations of the modified Bessel functions.

Hecke's result gives an Mellin--Barnes integral for the 
Bessel function of the third kind and also ties up with 
the modified Bessel functions. The most comprehensive work on the theory of Bessel functions is  
 \cite{watson} but a substantial account is available in  the book  of Lebedev \cite{lebedev} 
where on pages 117-118 we find a careful specification of the domain of validity for the Hankel functions. 

\paragraph*{Theorem 4.1 (Hecke):} For positive real numbers $p, q$ and arbitrary real $a$ and $b$, 
\begin{equation}\label{hecke}
\frac{1}{2\pi}\int_{-\infty}^{\infty}
\frac{\Gamma(a+c+it)\Gamma(b+c+it)\;dt}{p^{a+c+it}q^{b+c+it}} = 
\int_0^{\infty}\exp\Big(
-\frac{p}{s} - qs
\Big)  s^{b-a-1} ds               
\end{equation}
where $c$ is any real number such that $c + a > 0$ and $c + b > 0$.
\paragraph*{Proof:} Writing out $\Gamma(a+c+it)$ and $\Gamma(b+c+it)$ as integrals we get 
$$
\int_{-\infty}^{\infty}\frac{\Gamma(a+c+it)\Gamma(b+c+it)\;dt}{p^{a+c+it}q^{b+c+it}} = 
\int_0^{\infty}\int_0^{\infty} \frac{e^{(-u-v)}u^{a+c-1}v^{b+c-1}dudv}{p^{a+c}q^{b+c}}
\int_{-\infty}^{\infty}\exp(it\rho(u, v))\;dt
$$ 
where the phase function $\rho(u, v)$ is given by 
$\displaystyle{
\rho(u, v) = \log(uv) - \log(pq)
}$ 
and $\rho^{*}(\delta_0)$ along the hyperbola $uv = pq$ is given by the weighted measure $pqu^{-1}du$. Hence 
$$
\int_{-\infty}^{\infty}\frac{\Gamma(a+c+it)\Gamma(b+c+it)\;dt}{p^{a+c+it}q^{b+c+it}} = 
2\pi p^{b-a}\int_0^{\infty}\exp\Big(-u - \frac{pq}{u}\Big) u^{a-b-1} du
$$
The change of variables $u = ps$ now gives the desired result. 
\paragraph*{Remark:} We note that the special case $a = b = 0$ and $p = q = \sqrt{x}$ appears in \cite{voronoi}, p. 245 
in the following form:
$$
\frac{1}{2\pi i}\int_{c-i\infty}^{c+i\infty} \frac{\Gamma^2(s) ds}{x^s} = 
2\int_1^{\infty}\frac{e^{-2u\sqrt{x}} du}{\sqrt{u^2-1}}, \quad x > 0.
$$
The substitution $u + \sqrt{u^2-1} = t$  brings it in the form \eqref{hecke}.  

Taking $p = q$ we get the Mellin--Barnes integral for the Hankel functions $H_{p}^{(1)}(z)$ 
(see \cite{watson}, p. 180 or \cite{lebedev}, p. 108 and p. 118):  
\paragraph*{Theorem 4.2:} For the Bessel functions of the third kind we have the following Mellin--Barnes integral representation
\begin{equation}\label{hankel}
i\pi e^{i\lambda\pi/2}H_{\lambda}^{(1)}(2ix) = 
\frac{1}{2\pi}\int_{-\infty}^{\infty}\frac{\Gamma(\lambda+r+it)\Gamma(r+it)\;dt}{x^{\lambda+2r+2it}}, 
\quad x > 0 
\end{equation}
where $r > 0$ is arbitrary. 
\paragraph*{Proof:}  Setting $p = q = x$ and  $b-a = \lambda$ in \eqref{hecke} and denoting $a+c$ by $r$ we get, 
$$
\frac{1}{2\pi}\int_{-\infty}^{\infty}\frac{\Gamma(r+\lambda+it)\Gamma(r+it)\;dt}{x^{\lambda+2r+2it}} = \int_0^{\infty}\exp\Big(
-xt - \frac{x}{t}
\Big)t^{\lambda-1}dt
$$
The change of variables $t = e^v$ in the right hand side gives (see \cite{lebedev}, p. 118):
\begin{eqnarray*}
\frac{1}{2\pi}\int_{-\infty}^{\infty}
\frac{\Gamma(\lambda+r+it)\Gamma(r+it)\;dt}{x^{\lambda+2r+2it}} = 
2\int_0^{\infty}\exp(-2x\cosh u)\cosh(\lambda u)\; du 
= i\pi e^{i\lambda\pi/2}H_{\lambda}^{(1)}(2ix)
\end{eqnarray*}
A further change of variables gives the form in which the result appears in \cite{hecke}, p. 378.
\paragraph*{Remark:} The integral formula of Hecke continues to hold when $p$ and $q$ are complex  with positive real parts. 
We now derive the 
Mellin--Barnes integral representation for Bessel functions of the first kind (see \cite{bateman}, Volume - II, p. 21). 
\paragraph*{Theorem 4.3:} If $s$ is real positive, we have the following 
\begin{equation}\label{bessel}
J_{\lambda}(2s) = 
\frac{1}{2\pi}\int_{-\infty}^{\infty} \frac{\Gamma(\lambda+r+it)}{\Gamma(1-r-it)}s^{-\lambda-2r-2it}\;dt,  
\end{equation}
where $1 > r > -\mbox{Re}\;\lambda$ is arbitrary. 
\paragraph*{Proof:} 
In Hecke's formula we set 
$p = q = x + iy$ with $x > 0$ and $y = -s < 0$. We get in  the 
limit as $x\rightarrow 0+$  
\begin{equation}\label{hankel1}
\frac{1}{2\pi}\int_{-\infty}^{\infty} 
\frac{\Gamma(r+it)\Gamma(r+\lambda+it)}{|y|^{\lambda+2r+2it}}\exp(-\pi t + i\pi r)\;dt = i\pi H^{(1)}_\lambda(2|y|).
\end{equation}
Taking the complex conjugate of \eqref{hankel1} and performing the change of variables $t\mapsto -t$ we get 
(see \cite{lebedev}, p. 118):
\begin{equation}\label{hankel2}
\frac{1}{2\pi}\int_{-\infty}^{\infty}
\frac{\Gamma(r+it)\Gamma(r+\lambda+it)}{|y|^{\lambda+2r+2it}}\exp(\pi t - i\pi r)\;dt = -i\pi H^{(2)}_\lambda(2|y|). 
\end{equation}
Now subtracting \eqref{hankel2} from \eqref{hankel1} we get  
\begin{equation}\label{bessel1}
\frac{1}{2\pi}\int_{-\infty}^{\infty} 
\frac{\Gamma(r+it)\Gamma(r+\lambda+it)}{|y|^{\lambda+2r+2it}}\sin(r\pi + i\pi t) \;dt = 
\pi J_{\lambda}(2|y|).
\end{equation}
The desired result now follows on appealing to Euler's reflection formula.

An integral similar to Hecke's  with complex values of $p$ and $q$ is the following expression for Bessel functions 
of the first kind  due to Sonine (see \cite{watson}, p. 177):
\paragraph*{Theorem 4.4 (Sonine):} For $x > 0$ and Re$\;a > 0$ we have 
\begin{equation}\label{sonine}
J_{a}(x) = \frac{(x/2)^a}{2\pi}\int_{-\infty}^{\infty}(c+it)^{-a-1}\exp\Big(c+it - \frac{x^2}{4(c+it)}\Big)\;dt,  
\end{equation}
where $c > 0$ is arbitrary. 
\paragraph*{Proof:} We prove the result assuming that Re$\;a > 1/2$ and then appeal to analytic continuation.  
We take advantage of the formula for the Fourier transform of the Gaussian and  write the term 
$\exp(-x^2/4(c+it))$ in the form 
$$
\exp(-x^2/4(c+it)) = \Big(\frac{\pi}{c+it}\Big)^{-1/2}\int_{-\infty}^{\infty} \exp(-(c+it)y^2- ixy) \;dy 
$$
and the integral is then
$$
\frac{(x/2)^a}{2\Gamma(a+\frac{1}{2})\pi^{3/2}}
\int_{-\infty}^{\infty}\Gamma(a+\frac{1}{2})(c+it)^{-a-\frac{1}{2}}dt\int_{-\infty}^{\infty} 
\exp\Big((c+it)(1-y^2)-ixy\Big)\;dy
$$
Using the formula 
\begin{equation}\label{laplace_transform}
s^{-\lambda}\Gamma(\lambda) = \int_0^{\infty} e^{-su}u^{\lambda-1}\;du 
\end{equation}
the  integral may be rewritten as: 
$$
\frac{(x/2)^a}{2\Gamma(a+\frac{1}{2})\pi^{3/2}}\int_0^{\infty}du\int_{-\infty}^{\infty} u^{a-\frac{1}{2}}\exp(c(1-u-y^2)-ixy)\;dy
\int_{-\infty}^{\infty} \exp(it(1-u-y^2))\;dt.
$$
The integral with respect to $t$ is an oscillatory integral with phase $\rho(u, y) = u-1+y^2$ and we get 
using \eqref{funda-formula} 
$$
\frac{(x/2)^a}{\sqrt{\pi}\Gamma(a+\frac{1}{2})}
\langle
\rho^{*}(\delta_0), u^{a-\frac{1}{2}}\exp(c-cu - cy^2 - ixy)
\rangle;
$$
computed along $u = 1-y^2$. 
Note that since $u \geq 0$ only the part of the parabolic locus  $\rho(u, y) = 0$ with $y \in [-1, 1]$ 
is relevant. So we get finally 
$$
\frac{(x/2)^a}{\sqrt{\pi}\Gamma(a+\frac{1}{2})} \int_{-1}^1 (1-y^2)^{a-\frac{1}{2}}\cos xy dy,
$$ 
which is the integral form of Bessel function well suited for obtaining its asymptotic behavior 
(see \cite{hormander}, p. 286 or \cite{lebedev}, p. 114).

Appealing to analytic continuation, the formula of Sonine holds for all complex values of $z$ not on the negative real axis. 
We now examine the relationship between the formulas of Sonine and Hecke, namely from Sonine's formula we derive a 
representation of the modified Bessel function as a Hecke type integral \eqref{hecke}. Indeed, 
Sonine's formula immediately leads to the following formula for the modified 
Bessel functions (see \cite{lebedev}, p. 109 or \cite{rademacher}, p. 47):
\begin{equation}\label{basset}
x^{-p}I_p(2x) = \frac{1}{2\pi i}\int_{L} z^{-p-1}\exp\Big(
z + \frac{x^2}{z}
\Big)\;dz, \quad x > 0, 
\end{equation}
where $p> 0$, $c > 0$ and $L$ is the vertical line traced upwards and passing through $c$.  
From this we infer the following, 
where $K_p$ denotes the second modified Bessel function (see \cite{lebedev}, p. 108):  
\begin{eqnarray}
K_p(z) &=& \frac{\pi}{2}\frac{I_{-p}(z) - I_p(z)}{\sin p\pi},\quad p \notin \mathbb Z, \label{eqn: macdonald1}\\
K_n(z) &=& \lim_{p\rightarrow n}K_p(z),\quad n \in \mathbb Z \label{eqn: macdonald2}.
\end{eqnarray}
\paragraph*{Corollary 4.5:} For $p$ arbitrary real we have for all real positive $x$,  
\begin{eqnarray*}
I_{p}(2x) &=& \frac{-\sin p\pi}{\pi}\int_1^{\infty}\exp\Big(-xu - \frac{x}{u}
\Big)\frac{du}{u^{p+1}} + \frac{1}{\pi}\int_0^{\pi}\cos pt \exp(2x\cos t)\;dt, \\
K_p(2x) &=& \frac{1}{2}\int_0^{\infty} \exp\Big(
-ux -  \frac{x}{u}
\Big)\frac{du}{u^{p+1}}.\\
\end{eqnarray*}
\paragraph*{Proof:} Assume first that $p$ is real and $-1 < p < 1$. We can then deform the line $L$ in \eqref{basset} 
into a 
Hankel contour  consisting of the pair of parallel rays 
$L_{\pm} = \{z : z = t \pm i\epsilon, -\infty < t \leq -x\}$   
$(\epsilon > 0)$ and a circle of radius $x$ centered at the origin. 
The Hankel contour is traced so as to loop the origin in the counter-clockwise direction. 
In the limit as $\epsilon \rightarrow 0$ we get 
$$
I_p(2x) = -\frac{\sin p\pi}{\pi}\int_1^{\infty}\exp(-xu - \frac{x}{u})\frac{du}{u^{p+1}} + 
\frac{1}{2\pi}\int_{-\pi}^{\pi}\cos pt\exp(2x\cos t)\;dt,
$$
from which we get (by \eqref{eqn: macdonald1}) 
the stated integral representation for $K_p(x)$. The condition $-1 < p < 1$ can be relaxed 
by appealing to analytic continuation since the integrals are entire functions of $p$. The result clearly holds 
for complex values of  $x$ with positive real part. Via \eqref{hecke} we can immediately 
get a Mellin--Barnes integral for $K_p(x)$ (see \cite{rademacher}, p. 55). 
\paragraph*{Corollary 4.6:} For $x$ and $p$  real positive, we have the Mellin--Barnes representation 
$$
2x^pK_p(2x) = \frac{1}{2\pi i}\int_{L} \frac{\Gamma(z)\Gamma(z+p)\;dz}{x^{2z}}
$$
where $L$ is a vertical line in the right half-plane.

We now gather up a few other integrals whose proofs, by the method described here, reduce to complete triviality. 
\paragraph*{Theorem 4.7 (miscellaneous integrals):} The following formulas hold:
\begin{itemize} 
\item[(i)] Integral of Cahen: $\displaystyle{
\frac{1}{2\pi}\int_{-\infty}^{\infty} \Gamma(c + it)y^{-c-it}dt = e^{-y}, \quad 0 < c,\; 0 < y. 
}$ 
\item[(ii)]  Laplace's integral:  $\displaystyle{
\frac{1}{2\pi}\int_{-\infty}^{\infty} \exp(c + it)(c+it)^{-z}dt = (\Gamma(z))^{-1}
}$, $c > 0$, Re$\;z > 0$. 
\item[(iii)]  Fourier transform of $\Gamma(a+it)$: $\displaystyle{
\int_{-\infty}^{\infty} \Gamma(a+ it)e^{-it\xi} dt = 
2\pi\exp(-e^{\xi}+a\xi).
}$
\end{itemize}
\paragraph*{Remarks:} Laplace's integral appears in his famous work on probability theory (see \cite{laplace}, p. 137). 
By deforming the contour of Laplace's integral it is easy to deduce Hankel's integral 
representation for $\Gamma(z)^{-1}$ which is valid for all $z\in \mathbb C$ (see \cite{nielsen}, pp. 147-148). 
The Fourier transform (iii) appears in \cite{pincherle}, p. 235. Cahen's integral is valid for complex values of $y$ 
subject to the conditions $|\mbox{arg} (y)| < \pi/2$ and $y \neq 0$.

We have seen that the method used in this paper handles many integrals which belong to the class of  
Meijer G functions \cite{mainardi}, \cite{bateman} (volume 1). 
The precise subclass of Meijer G functions to which this method 
would be applicable is a matter for further investigation.  
\paragraph*{Acknowledgements:} The author thanks his teacher Professor Paul Garrett at the University of 
Minnesota for reading the manuscript 
and suggesting changes. 
Indeed to Professor Garrett, the author owes a great deal more than words can express.  
The author thanks Professor Steven Krantz for his encouraging remarks after 
 reading the manuscript and for providing the much needed moral boost.      
The author thanks the referees for pointing out some references and for suggesting changes that has resulted in a 
considerable improvement of the paper. 

\end{document}